\documentclass[11pt]{article}
\usepackage{amssymb}

\usepackage{amsmath}
\usepackage{graphicx}
\usepackage{theorem}
\usepackage{here}

\newtheorem{thm}{Theorem}[section]
\newtheorem{lem}{Lemma}[section]

\theorembodyfont{\rmfamily}

\makeatletter

\@addtoreset{equation}{section}
\makeatother
\def\al{{\alpha}}

\def\ga{{\gamma}}

\def\la{{\lambda}}
\def\si{{\sigma}}

\def\bde{{\text{\boldmath $\delta$}}}

\def\bom{{\text{\boldmath $\omega$}}}
\def\bth{{\text{\boldmath $\theta$}}}

\def\bmu{{\text{\boldmath $\mu$}}}
\def\bnu{{\text{\boldmath $\nu$}}}

\def\bnabla{{\text{\boldmath $\nabla$}}}

\def\sih{{\hat \si}}

\def\bthh{{\widehat \bth}}

\def\bmuh{{\widehat \bmu}}

\def\Si{{\Sigma}}

\def\bSi{{\text{\boldmath $\Si$}}}

\def\a{{\text{\boldmath $a$}}}

\def\h{{\text{\boldmath $h$}}}

\def\j{{\text{\boldmath $j$}}}

\def\x{{\text{\boldmath $x$}}}
\def\y{{\text{\boldmath $y$}}}

\def\A{{\text{\boldmath $A$}}}

\def\C{{\text{\boldmath $C$}}}

\def\H{{\text{\boldmath $H$}}}
\def\I{{\text{\boldmath $I$}}}

\def\Q{{\text{\boldmath $Q$}}}

\def\V{{\text{\boldmath $V$}}}
\def\W{{\text{\boldmath $W$}}}
\def\X{{\text{\boldmath $X$}}}
\def\Y{{\text{\boldmath $Y$}}}

\def\Nc{{\cal N}}

\def\Re{{\mathbb{R}}}
\def\tr{{\rm tr\,}}

\def\[{{\text{\boldmath $[$}}}
\def\]{{\text{\boldmath $]$}}}
\def\et{{\it et\, al.}}
\def\zero{{\bf\text{\boldmath $0$}}}

\def\|{{\,|\,}}

\def\/{{\Bigr/\!\!}}
\def\non{{\nonumber}}

\def\bnuh{{\widehat \bnu}}
\def\ddo{{\overline d}}

\begin{document}
\title{Bayes Minimax Competitors of Preliminary Test Estimators in $k$ Sample Problems}

\author{
Ryo Imai\footnote{Graduate School of Economics, University of Tokyo, 
7-3-1 Hongo, Bunkyo-ku, Tokyo 113-0033, JAPAN.\newline{
E-Mail: imronuronnar13@g.ecc.u-tokyo.ac.jp  }},
Tatsuya Kubokawa\footnote{Faculty of Economics, University of Tokyo, 
7-3-1 Hongo, Bunkyo-ku, Tokyo 113-0033, JAPAN. \newline{
E-Mail: tatsuya@e.u-tokyo.ac.jp}}
and
Malay Ghosh\footnote{Department of Statistics, University of Florida, 102 Griffin-Floyd Hall, Gainesville, Florida 32611.\newline{
E-Mail: ghoshm@stat.ufl.edu  }}
}
\maketitle
\begin{abstract}
In this paper, we consider the estimation of a mean vector of a multivariate normal population where the mean vector is suspected to be nearly equal to mean vectors of $k-1$ other populations.
As an alternative to the preliminary test estimator based on the test statistic for testing hypothesis of equal means, we derive empirical and hierarchical Bayes estimators which shrink the sample mean vector  toward a pooled mean estimator given under the hypothesis.
The minimaxity of those Bayesian estimators are shown, and their performances are investigated by simulation.

\par\vspace{4mm}
\end{abstract}

\section{Introduction}

Suppose that there are $k$ laboratories, say Laboratory L$_1$, $\ldots$, L$_k$, and a certain instrument is designed to measure several characteristics at each laboratory and several vector-valued measurements are recorded.
Also, suppose that we want to estimate the population mean of Laboratory $L_1$.
When similar instruments are used at $k$ laboratories, it is suspected that $k$ population means are nearly equal, in which case, the sample means of the other laboratories L$_2$, $\ldots$, L$_k$ are used to produce more efficient estimators than the sample mean based on data from only L$_1$.
This problem was studies by Ghosh and Sinha (1988) and recently revisited by Imai, Kubokawa and Ghosh (2017) in the framework of simultaneous estimation of $k$ population means.

\medskip
The $k$ sample problem described above is expressed in the following canonical form:
$p$-variate random vectors $\X_1, \ldots, \X_k $ and positive scalar random variable $S$ are mutually independently distributed as
\begin{equation}
\begin{split}
\X_i \sim& \Nc_p(\bmu_i, \si^2\V_i),\quad \text{for}\ i=1, \ldots, k,\\
S/\si^2\sim& \chi_n^2,
\end{split}
\label{eqn:model}
\end{equation}
where the $p$-variate means $\bmu_1, \ldots, \bmu_k$ and the scale parameter $\si^2$ are unknown, and $\V_1, \ldots, \V_k$ are $p\times p$ known and positive definite matrices.
In this model, we consider to estimate $\bmu_1$ relative to the quadratic loss function
\begin{equation}
L(\bde_1, \bmu_1, \si^2) = \Vert\bde_1-\bmu_1\Vert_{\Q}^2/\si^2
=(\bde_1-\bmu_1)^\top\Q (\bde_1-\bmu_1)/\si^2,
\label{eqn:loss}
\end{equation}
where $\Vert \a \Vert_\A^2=\a^\top \A \a$ for the transpose $\a^\top$ of $\a$, $\Q$ is a positive definite and known matrix, and $\bde_1$ is an estimator of $\bmu_1$.
Estimator $\bde_1$ is evaluated by the risk function $R(\bom, \bde_1)=E[L(\bde_1, \bmu_1, \si^2)]$ for $\bom=(\bmu_1, \ldots, \bmu_k, \si^2)$, a set of unknown parameters.

\medskip
In this paper, we consider the case that the means $\bmu_i$'s are suspected to be nearly equal, namely close to the hypothesis 
$$
H_0 : \bmu_1=\cdots = \bmu_k.
$$
A classical approach towards solution to this problem is the development of a preliminary test estimator which uses the pooled mean estimator upon acceptance of  the null hypothesis $H_0$, and uses separate mean estimators upon rejection of $H_0$.
The test statistic for $H_0$ is 
\begin{equation}
\label{eqn:test statistic}
F=\sum_{i=1}^k \Vert \X_i-\bnuh\Vert_{\V_i^{-1}}^2/S,
\end{equation}
where $\bnuh$ is the pooled estimator defined as
\begin{equation}
\bnuh=\A \sum_{i=1}^k \V_i^{-1}\X_i, \quad
\A=\Big(\sum_{i=1}^k\V_i^{-1}\Big)^{-1}.
\label{eqn:nuh}
\end{equation}
The preliminary test estimator of $\bmu_1$ is 
\begin{equation}
\bmuh_1^{PT}= 
\left\{\begin{array}{ll} \X_1 & \ \text{if}\ F>(p(k-1)/n)F_{p(k-1), n, \al}\\
\bnuh & \ \text{otherwise}
\end{array}\right.
\label{eqn:PT}
\end{equation}
where $F_{p(k-1), n, \al}$ is the upper $\al$ point of the F distribution with $(p(k-1), n)$ degrees of freedom, and $\bnuh$ is the pooled estimator given in (\ref{eqn:nuh}).
However, the preliminary test estimator is not smooth and does not necessarily improve on $\X_1$.

\medskip
As an alternative approach to the preliminary test estimator, we consider a Bayesian method under the prior distributions of $\bmu_i$ having a common mean $\bnu$ as a hyper-parameter.
Empirical and hierarchical Bayes estimators are derived when the uniform prior distribution is assumed for $\bnu$.
We also provide empirical Bayes estimator under the assumption of the normal prior distribution for $\bnu$.
It is shown that these Bayesian estimators improve on $\X_1$, namely, are minimax.

\medskip
The topic treated in this paper is related to the so-called Stein problem in simultaneous estimation of multivariate normal means.
This problem has long been of interest in the literature.
See, for example, Stein (1956, 1981), James and Stein (1961), Strawderman (1971, 1973), Efron and Morris (1973, 1976)  and Berger (1985).
For recent articles exrending to prediction and high-dimensional problems, see Komaki (2001), Brown, George and Xu (2008) and Tsukuma and Kubokawa (2015).
As articles related to this paper, Sclove, Morris and Radhakrishnan (1972) showed the inadmissibility of the preliminary test estimator, Smith (1973) provided Bayes estimators in one-way and two-way models,  Ghosh and Sinha (1988) derived hierarchical and empirical Bayes estimators with minimaxity, and Sun (1996) provided Bayesian minimax estimators in multivariate two-way random effects models.

\medskip
In Section \ref{sec:2}, we treat two classes of shrinkage estimators, say Class 1 and Class 2.
Estimators in Class 1 shrink $\X_1$ toward the pooled estimator $\bnuh$, and estimators in Class 2 incorporate a part of shrinking $\bnuh$ in addition to estimators in Class 1.
Ghosh and Sinha (1988) obtained conditions for minimaxity of estimators in Class 1 in the two sample problem.
For $k=2$, the test statistic $F$ can be simply written as $F=(\X_1-\X_2)^\top(\V_1+\V_2)^{-1}(\X_1-\X_2)/S$, while $F$ is more complicated as $F=\sum_{i=1}^k (\X_i-\bnuh)^\top\V_i^{-1}(\X_i-\bnuh)/S$ in the $k$ sample case.
This means that the minimaxity is harder to establish in the $k$ sample problem.
The key tool in the extension is the inequality given in Lemma \ref{lem:in}.
Using the inequality, we obtain conditions for minimaxity of shrinkage estimators in Class 1 and Class 2. 

\medskip
In Section \ref{sec:est}, we suggest three kinds of Bayesian methods for estimation of $\bmu_1$.
Empirical and hierarchical Bayes estimators are derived under the assumption of the uniform prior distribution for $\bnu$.
Because these estimators belong to Class 1, we can provide conditions for their minimaxity using the result in Section \ref{sec:2}.
Also empirical Bayes estimator is obtained under the normal prior distribution for $\bnu$ and belongs to Class 2.
Conditions for the minimaxity are given from the result in Section \ref{sec:2}. 
The performances of those Bayesian procedures are investigated by simulation in Section \ref{sec:sim}. 
The results in Section \ref{sec:2} are extended in Section \ref{sec:4} to the problem of estimating a linear combination of $\bmu_1, \ldots, \bmu_k$ under the quadratic loss. 
Concluding remarks are given in Section \ref{sec:remark}.

\section{Classes of Minimax Estimators}
\label{sec:2}

Motivated from the preliminary test estimator (\ref{eqn:PT}), we consider a class of estimators shrinking $\X_1$ toward the pooled estimator $\bnuh$, given by
\begin{equation}
\bmuh_1(\phi) = \X_1 - {\phi(F, S)\over F}(\X_1-\bnuh),
\label{eqn:class1}
\end{equation}
where $\phi(F,S)$ is an absolutely continuous function, $F$ is the test statistic given in (\ref{eqn:test statistic}) and $\bnuh$ is the pooled estimator $\bnuh$  given in (\ref{eqn:nuh}).

We can derive conditions on $\phi(F,S)$ for minimaxity of $\bmuh_1(\phi)$ in Theorem \ref{thm:1}, which is proved in the end of this section.

\begin{thm}
\label{thm:1}
Assume the conditions
\begin{equation}
{\rm Ch}_{\rm max}\{(\V_1-\A)\Q\}\neq0\quad {\rm and}\quad
\tr\{(\V_1-\A)\Q\}/{\rm Ch}_{\rm max}\{(\V_1-\A)\Q\}>2.
\label{eqn:21a}
\end{equation}
Then, the estimator $\bmuh_1(\phi)$ is minimax relative  to the quadratic loss $(\ref{eqn:loss})$ if $\phi(F,S)$ satisfies the following conditions:

{\rm (a)}\ $\phi(F,S)$ is non-decreasing in $F$ and non-increasing in $S$.

{\rm (b)}\ $0<\phi(F,S)\leq 2 [\tr\{(\V_1-\A)\Q\}/{\rm Ch}_{\rm max}\{(\V_1-\A)\Q\}-2]/(n+2)$, where ${\rm Ch}_{\rm max}(\C)$ denotes the maximum characteristic value of matrix $\C$.
\end{thm}

This theorem provides an extension of Ghosh and Sinha (1988) to the $k$-sample problem.
In the case of $\V_1=\cdots=\V_k=\Q^{-1}$, the conditions in (\ref{eqn:21a}) are expressed as $k\not= 1$ and $p>2$, and the latter condition is well known in the Stein problem.

We next consider the class of double shrinkage estimators
\begin{equation}
\bmuh_1(\phi, \psi) = \X_1 - {\phi(F, S)\over F}(\X_1-\bnuh)- {\psi(G, S)\over G}\bnuh,
\label{eqn:class2}
\end{equation}
where $\phi(F,S)$ and $\psi(G,S)$ are absolutely continuous functions, and 
$$
G=\Vert\bnuh\Vert_{\A^{-1}}^2/S.
$$

\begin{thm}
\label{thm:2}
Assume condition $(\ref{eqn:21a})$ and 
\begin{equation}
{\rm Ch}_{\rm max}(\A\Q)\neq0 \quad {\rm and}\quad \tr(\A\Q)/{\rm Ch}_{\rm max}(\A\Q)>2.
\label{eqn:22a}
\end{equation}
Then, the double shrinkage estimator $\bmuh_1(\phi, \psi)$ in $(\ref{eqn:class2})$ is minimax relative  to the quadratic loss $(\ref{eqn:loss})$ if $\phi(F,S)$ and $\psi(G,S)$ satisfy the following conditions:

{\rm (a)}\ $\phi(F,S)$ is non-decreasing in $F$ and non-increasing in $S$.

{\rm (b)}\ $0<\phi(F,S)\leq [\tr\{(\V_1-\A)\Q\}/{\rm Ch}_{\rm max}\{(\V_1-\A)\Q\}-2]/(n+2)$.

{\rm (c)}\ $\psi(G,S)$ is non-decreasing in $G$ and non-increasing in $S$.

{\rm (d)}\ $0<\psi(G,S)\leq \{\tr(\A\Q)/{\rm Ch}_{\rm max}(\A\Q)-2\}/(n+2)$.
\end{thm}

For the proofs, the Stein identity due to Stein (1981) and the chi-square identity due to Efron and Morris (1976) are useful.
See also Bilodeau and Kariya (1989) for a multivariate version of the Stein identity.

\begin{lem}
\label{lem:identity}
{\rm (1)} \ Assume that $\Y=(Y_1, \ldots, Y_p)^\top$ is a $p$-variate random vector having $\Nc_p(\bmu, \bSi)$ and that $\h(\cdot)$ is an absolutely continuous function from $\Re^p$ to $\Re^p$.
Then, the Stein identity is given by
\begin{equation}
E[(\Y-\bmu)^\top \h(\Y)]=E\big[ \tr\big\{ \bSi \bnabla_\Y \h(\Y)^\top\big\}\big],
\label{eqn:stein}
\end{equation}
provided the expectations in both sides exist, where $\bnabla_\Y=(\partial /\partial Y_1, \ldots, \partial /\partial Y_p)^\top$.

\medskip
{\rm (2)}\ Assume that $S$ is a random variable such that $S/\si^2\sim \chi_n^2$ and that $g(\cdot)$ is an absolutely continuous function from $\Re$ to $\Re$.
Then, the chi-square identity is given by
\begin{equation}
E[S g(S)] = \si^2 E[n g(S) + 2Sg'(S)],
\label{eqn:chi-square}
\end{equation}
provided the expectations in both sides exist.
\end{lem}

{\bf Proof of Theorem \ref{thm:1}}\ \ 
The risk function is decomposed as
\begin{align}
R(\bom, \bmuh_1(\phi))
=&
E[\Vert \X_1-\bmu_1\Vert_\Q^2/\si^2]
-2 E\Big[(\X_1-\bmu_1)^\top\Q(\X_1-\bnuh){\phi\over \si^2 F}\Big]
 + {1\over \si^2} E\Big[ {\phi^2\over F^2}\Vert \X_1-\bnuh\Vert_\Q^2\Big]
\non\\
=& I_1 - 2 I_2 + I_3, \quad\quad \text{(say)}
\label{eqn:p1}
\end{align}
for $\phi=\phi(F,S)$ and $\psi=\psi(G,S)$.

\medskip
It is easy to see $I_1=\tr(\V_1\Q)$.
Note that $F$ is expressed as $F=(\sum_{i=1}^k\X_i^\top\V_i^{-1}\X_i-\bnuh^\top\A^{-1}\bnuh)/S$.
Letting $\bnabla_1=\partial/\partial \X_1$, we have $\bnabla_1 F=2\V_1^{-1}(\X_1-\bnuh)/S$ and
\begin{align*}
	\bnabla_1\Big\{(\X_1-\bnuh)^\top {\phi(F,S)\over F} \Big\}=(\I-\V_1^{-1}\A) {\phi\over F}+2\V_1^{-1}(\X_1-\bnuh)(\X_1-\bnuh)^\top\Big\{ - {\phi\over F^2} +{ \phi_F\over F}\Big\}\frac{1}{S},
\end{align*}
where $\phi_F=(\partial /\partial F)\phi(F,S)$. Then from (\ref{eqn:stein}), it is seen that
\begin{align}
I_2 =& E\Big[\tr\Big[\Q\V_1\bnabla_1 \Big\{(\X_1-\bnuh)^\top {\phi(F,S)\over F}\Big\}\Big]\Big]
\non\\
=&
E\Big[ \tr\{ \Q\V_1 (\I-\V_1^{-1}\A) \} {\phi\over F} + 2{\tr\{\Q\V_1\V_1^{-1}(\X_1-\bnuh)(\X_1-\bnuh)^\top\}\over S}\Big\{ - {\phi\over F^2} +{ \phi_F\over F}\Big\} \Big]
\non\\
=&
E\Big[ \tr\{ (\V_1-\A)\Q\}{\phi \over F} - 2 B{\phi \over F} + 2 B\phi_F\Big], 
\label{eqn:I2}
\end{align}
where
\begin{equation}
B= (\X_1-\bnuh)^\top\Q(\X_1-\bnuh)/\sum_{j=1}^k (\X_j-\bnuh)^\top\V_j^{-1}(\X_j-\bnuh).
\label{eqn:B}
\end{equation}
Also, from (\ref{eqn:chi-square}), it is observed that
\begin{align}
I_3=&
{1\over \si^2} E\Big[{S\over F}B\phi^2\Big] 
\non\\
=&  E\Big[{n\over F}B\phi^2+ 2SB\Big(-{F\over S}\Big) \Big\{ - {\phi^2\over F^2}+2{\phi\phi_F\over F}\Big\} + 2S B{\phi\phi_S\over F} \Big]
\non\\
=&
E\Big[B{n+2\over F}\phi^2- 4B\phi\phi_F + 4S B{\phi\phi_S\over F} \Big].
\label{eqn:I3}
\end{align} 
Thus, the risk function is expressed as 
\begin{align}
R(\bom, \bmuh_1(\phi))
=&
E\Big[ \tr(\V_1\Q)  + {\phi\over F}[(n+2)B\phi - 2 \tr\{(\V_1-\A)\Q\}+4B]
\non\\
&
-4B\phi_F-4B\phi\phi_F+4S B{\phi\phi_S\over F} \Big].
\label{eqn:p2}
\end{align}
Using Lemma \ref{lem:in} given below, we can see that
$$
B \leq {(\X_1-\bnuh)^\top\Q(\X_1-\bnuh) \over (\X_1-\bnuh)^\top(\V_1-\A)^{-1}(\X_1-\bnuh)}\leq {\rm Ch}_{\rm max}((\V_1-\A)\Q),
$$
which implies that 
\begin{align}
R(&\bom, \bmuh_1(\phi)) - \tr(\V_1\Q)
\non\\
\leq&
E\Big[{\rm Ch}_{\rm max}((\V_1-\A)\Q){\phi\over F}\Big\{(n+2)\phi - 2 {\tr\{(\V_1-\A)\Q\}\over {\rm Ch}_{\rm max}((\V_1-\A)\Q)}+4\Big\}
\non\\
&-4B\phi_F-4B\phi\phi_F+4S B{\phi\phi_S\over F} \Big].
\label{eqn:p3}
\end{align}
Hence, $R(\bom, \bmuh_1(\phi)) \leq \tr(\V_1\Q)$ under the conditions in Theorem \ref{thm:1}.
\hfill$\Box$

\bigskip
\begin{lem}
\label{lem:in}
It holds that
\begin{equation}
\sum_{j=1}^k(\x_j-\bnuh)^\top\V_j^{-1}(\x_j-\bnuh) \geq (\x_1-\bnuh)^\top(\V_1-\A)^{-1}(\x_1-\bnuh).
\label{eqn:in}
\end{equation}
\end{lem}

{\bf Proof}.\  \ 
Let $\C=(\V_1-\A)^{-1}$ and $\x_*=\sum_{j=2}^k\V_j^{-1}\x_j$.
Then, $\bnuh=\A\V_1^{-1}\x_1+\A\x_*$ and $\x_1-\bnuh=\C^{-1}\V_1^{-1}\x_1-\A\x_*$.
The RHS of (\ref{eqn:in}) is rewritten as
\begin{equation}
(\x_1-\bnuh)^\top\C(\x_1-\bnuh)=
\x_1^\top\V_1^{-1}\C^{-1}\V_1^{-1}\x_1
- 2\x_1^\top\V_1^{-1}\A\x_* + \x_*^\top\A\C\A\x_*.
\label{eqn:in1}
\end{equation}
On the other hand, it can be observed that
\begin{align*}
(\x_1-\bnuh)^\top\V_1^{-1}(\x_1-\bnuh)=&
\x_1^\top\V_1^{-1}\C^{-1}\V_1^{-1}\C^{-1}\V_1^{-1}\x_1
- 2\x_1^\top\V_1^{-1}\C^{-1}\V_1^{-1}\A\x_* + \x_*^\top\A\V_1^{-1}\A\x_*,
\\
\sum_{j=2}^k(\x_j-\bnuh)^\top\V_j^{-1}(\x_j-\bnuh)
=&
\sum_{j=2}^k\x_j^\top\V_j^{-1}\x_j
-2\x_*^\top\A\V_1^{-1}\x_1 - 2\x_*^\top\A\x_*
\\
&
+\x_1^\top\V_1^{-1}\A(\A^{-1}-\V_1^{-1})\A\V_1^{-1}\x_1
+ 2\x_1\V_1^{-1}\A(\A^{-1}-\V_1^{-1})\A\x_* 
\\
&+ \x_*^\top\A(\A^{-1}-\V_1^{-1})\A\x_*,
\end{align*}
which gives
\begin{align}
\sum_{j=1}^k&(\x_j-\bnuh)^\top\V_j^{-1}(\x_j-\bnuh) 
\non\\
=&
\x_1^\top\V_1^{-1}\{\C^{-1}\V_1^{-1}\C^{-1}+\A(\A^{-1}-\V_1^{-1})\A\}\V_1^{-1}\x_1
\non\\
&-2\x_1^\top\V_1^{-1}\{\C^{-1}\V_1^{-1}+\I-\A(\A^{-1}-\V_1^{-1})\}\A\x_*
\non\\
&-\x_*^\top\A\x_* + \sum_{j=2}^k\x_j^\top\V_j^{-1}\x_j.
\label{eqn:in2}
\end{align}
It is noted that $\C^{-1}\V_1^{-1}\C^{-1}+\A(\A^{-1}-\V_1^{-1})\A=\V_1-\A$, $\C^{-1}\V_1^{-1}+\I-\A(|A^{-1}-\V_1^{-1})=\I$ and $\A\C\A+\A=\A(\V_1-\A)^{-1}\V_1=(\A^{-1}-\V_1^{-1})^{-1}=(\sum_{j=2}^k\V_j^{-1})^{-1}$.
From (\ref{eqn:in1}) and (\ref{eqn:in2}), the inequality (\ref{eqn:in}) is equivalent to
\begin{equation}
\sum_{j=2}^k \x_j^\top\V_j^{-1}\x_j \geq \Big(\sum_{j=2}^k\V_j^{-1}\x_j\Big)^\top \Big(\sum_{j=2}^k\V_j^{-1}\Big)^{-1} \Big(\sum_{j=2}^k\V_j^{-1}\x_j\Big).
\label{eqn:in3}
\end{equation}
To show the inequality (\ref{eqn:in3}), let $\bnu_*=(\sum_{j=2}^k\V_j^{-1})^{-1} (\sum_{j=2}^k\V_j^{-1}\x_j)$.
Then, it can be seen that
\begin{align*}
\sum_{j=2}^k& \x_j^\top\V_j^{-1}\x_j - \Big(\sum_{j=2}^k\V_j^{-1}\x_j\Big)^\top \Big(\sum_{j=2}^k\V_j^{-1}\Big)^{-1} \Big(\sum_{j=2}^k\V_j^{-1}\x_j\Big)
\\
=& \sum_{j=2}^k \x_j^\top\V_j^{-1}\x_j - \bnu_*^\top \Big(\sum_{j=2}^k\V_j^{-1}\Big) \bnu_*
\\
=& \sum_{j=2}^k (\x_j-\bnu_*)^\top\V_j^{-1}(\x_j-\bnu_*),
\end{align*}
which is nonnegative, and the proof of Lemma \ref{lem:in} is complete.
\hfill$\Box$

\bigskip
{\bf Proof of Theorem \ref{thm:2}}\ \ 
The risk function of $\bmuh_1(\phi,\psi)$ is 
\begin{align*}
R(\bom, \bmuh_1(\phi,\psi))=&
R(\bom,\bmuh_1(\phi)) 
-2 {1\over \si^2}E\Big[ (\X_1-\bmu_1)^\top\Q\bnuh{\psi\over G}\Big]
\non\\
&
+2{1\over \si^2}E\Big[ {\phi \psi \over FG}(\X_1-\bnuh)^\top\Q\bnu\Big]
+ E\Big[ {\psi^2\over \si^2G^2}\bnuh^\top\Q\bnuh\Big],
\end{align*}
for $\phi=\phi(F,S)$ and $\psi=\psi(G,S)$.
Because of $(\X_1-\bnuh)^\top\Q\bnuh \leq \{(\X_1-\bnuh)^\top\Q(\X_1-\bnuh)\}^{1/2} \{\bnuh^\top\Q\bnuh\}^{1/2}$, we have
\begin{align*}
{2\over \si^2}{\phi \psi \over FG}(\X_1-\bnuh)^\top\Q\bnuh
\leq& {2S^2\over \si^2} \Big\{\phi {\{(\X_1-\bnuh)^\top\Q(\X_1-\bnuh)\}^{1/2}\over \sum_{j=1}^k\Vert\X_j-\bnuh\Vert_{\V_j^{-1}}^2}\Big\}  \Big\{ \psi {\{\bnuh^\top\Q\bnuh\}^{1/2}\over \bnuh^\top\A^{-1}\bnuh}\Big\}
\\
\leq&
{2S^2\over \si^2} \Big\{\phi^2 {(\X_1-\bnuh)^\top\Q(\X_1-\bnuh)\over 2\{\sum_{j=1}^k\Vert\X_j-\bnuh\Vert_{\V_j^{-1}}^2\}^2 }+ \psi^2 {\bnuh^\top\Q\bnuh\over 2(\bnuh^\top\A^{-1}\bnuh)^2}\Big\}
\\
=& {S\over \si^2}{\phi^2\over F} B+ {S\over \si^2}{\psi^2\over G} {\bnuh^\top\Q\bnuh\over \bnuh^\top\A^{-1}\bnuh},
\end{align*}
for $B$ defined in (\ref{eqn:B}).
Thus, 
\begin{align}
R(\bom, \bmuh_1(\phi,\psi))=&
R(\bom,\bmuh_1(\phi)) 
-2 {1\over \si^2}E\Big[ (\X_1-\bmu_1)^\top\Q\bnuh{\psi\over G}\Big]
\non\\
&
+ 2E\Big[ {S\over \si^2}{\psi^2\over G} W\Big]
+E\Big[ {S\over \si^2}{\phi^2\over F} B\Big]
\non\\
=& R(\bom,\bmuh_1(\phi)) -2J_1 + 2J_2 +J_3, \quad {\rm (say)}
\label{eqn:p5}
\end{align}
where $W={\bnuh^\top\Q\bnuh/ \bnuh^\top\A^{-1}\bnuh}$.
The Stein identity is applied to rewrite $J_1$ as
\begin{equation}
J_1=E\Big[ {\psi\over G}\tr(\A\Q) - 2 W{\psi\over G}+ 2 W\psi_G \Big].
\label{eqn:J1}
\end{equation}
The chi-square identity is used to rewrite $J_2$ as
\begin{equation}
J_2=
E\Big[ (n+2)W{\psi^2\over G} -4W\psi\psi_G + 4WS{\psi\psi_S\over G}\Big].
\label{eqn:J2}
\end{equation}
The calculation of $J_3$ is given in (\ref{eqn:I3}), and from (\ref{eqn:p2}), (\ref{eqn:p5}), (\ref{eqn:J1}) and (\ref{eqn:J2}), it follows that
\begin{align}
R(\bom&, \bmuh_1(\phi,\psi))-\tr(\V_1\Q) 
\non\\
=&
E\Big[  {\phi\over F}\big\{2(n+2)B\phi - 2 \tr\{(\V_1-\A)\Q\}+4B\big\}
+ {\psi\over G}\big\{ 2(n+2)W\psi - 2\tr(\A\Q) +4W \big\}
\non\\
&
-4B\phi_F-8B\phi\phi_F+8BS{\phi\phi_S\over F}
-4 W\psi_G  -8W\psi\psi_G + 8WS{\psi\psi_S\over G}
\Big].
\label{eqn:p6}
\end{align}
Because $W\leq {\rm Ch}_{\rm max}(\A\Q)$, it can be verified that $R(\bom, \bmuh_1(\phi,\psi))\leq\tr(\V_1\Q) $ under the conditions (a)-(d) in Theorem \ref{thm:2}, which is proved.
\hfill$\Box$

\section{Hierarchical and Empirical Bayes Minimax Estimators}
\label{sec:est}

\subsection{ Empirical Bayes estimator under the uniform prior for $\bnu$}
\label{sec:3.1}

We begin with assuming the prior distribution 
\begin{equation}
\begin{split}
\bmu_i\mid \bnu, \tau^2 \sim& \Nc_p(\bnu, \tau^2\V_i), \quad \text{for}\ i=1, \ldots, k,
\\
\bnu \sim& \text{Uniform}(\Re^p),
\end{split}
\label{eqn:prior1}
\end{equation}
where Uniform$(\Re^p)$ denotes the improper uniform distribution over $\Re^p$, and $\tau^2$ is an unknown parameter.
The posterior distribution of $\bmu_i$ given $\X_i$ and $\bnu$, the posterior distribution of $\bnu$ given $\X_1, \ldots, \X_k$ and the marginal distribution of $\X_1, \ldots, \X_k$ are
\begin{equation}
\begin{split}
\bmu_i \mid \X_i, \bnu, \tau^2, \si^2 \sim& \Nc_p\Big( \bmuh_i^*(\si^2, \tau^2,\bnu), (\si^{-2}+\tau^{-2})^{-1}\Big),\quad  i=1, \ldots, k,
\\
\bnu \mid \X_1, \ldots, \X_k, \tau^2, \si^2 \sim& \Nc_p( \bnuh, (\tau^2+\si^2)\A),
\\
f_\pi (\x_1, \ldots, \x_k \mid \tau^2, \si^2) \propto& {1\over (\tau^2+\si^2)^{p(k-1)/2}}\exp\Big\{ -{\sum_{i=1}^k\Vert\x_i-\bnuh\Vert_{\V_i^{-1}}^2 \over 2(\tau^2+\si^2)} \Big\},
\end{split}
\label{eqn:posterior1}
\end{equation}
where $\bmuh_i^*(\si^2, \tau^2,\bnu) = \X_i - \{\si^2/(\tau^2+\si^2)\}(\X_i-\bnu)$.
Then, the Bayes estimator of $\bmu_1$ is 
\begin{equation}
\bmuh_1^B(\si^2, \tau^2)=E[\bmuh_1^*(\si^2, \tau^2,\bnu) \mid \X_1, \ldots,\X_k]
=\X_1 - {\si^2\over \tau^2+\si^2}(\X_1-\bnuh).
\label{eqn:Bayes2}
\end{equation}

Because $\tau^2+\si^2$ and $\si^2$ are unknown, we estimate $\tau^2+\si^2$ by $\sum_{i=1}^k \Vert \X_i-\bnuh\Vert_{\V_i^{-1}}^2/\{p(k-1)-2\}$ from the marginal likelihood in (\ref{eqn:posterior1}).
When $\si^2$ is estimated by $\sih^2=S/(n+2)$, the resulting empirical Bayes estimator is
\begin{equation}
\bmuh_1^{EB} = \X_1 - \min\Big( {a_0 \over F}, 1\Big)(\X_1-\bnuh),
\label{eqn:EB1}
\end{equation}
for $a_0=\{p(k-1)-2\}/(n+2)$.
It follows from Theorem \ref{thm:1} that the empirical Bayes estimator $\bmuh_1^{EB}$ is minimax for $0<a_0 \leq 2 [\tr\{(\V_1-\A)\Q\}/{\rm Ch}_{\rm max}\{(\V_1-\A)\Q\}-2]/(n+2)$.

\subsection{Hierarchical Bayes minimax estimator under the uniform prior for $\bnu$}
\label{sec:3.2}

We consider the prior distribution for $\si^2$ in (\ref{eqn:model}) and $\tau^2$ in (\ref{eqn:prior1}), namely, in addition of (\ref{eqn:prior1}), we assume that
\begin{equation}
\begin{split}
\pi(\tau^2 \mid \si^2) \propto& \Big({\si^2\over \tau^2+\si^2}\Big)^{a +1},\\
\pi(\si^2) \propto& (\si^2)^{c -2}, \quad \text{for}\ \si^2\leq 1/L,
\end{split}
\label{eqn:prior2}
\end{equation}
where $a $ and $c $ are constants, and $L$ is a positive constant.
From (\ref{eqn:posterior1}), the posterior distribution of $(\tau^2, \si^2)$ given $\X_1, \ldots, \X_k, S$ is
\begin{align*}
\pi&(\tau^2, \si^2\mid \x_1, \ldots, \x_k, S)\\
&
\propto  \Big({\si^2\over \tau^2+\si^2}\Big)^{p(k-1)/2+a +1}\Big({1\over \si^2}\Big)^{\{n+p(k-1)\}/2+2-c}\exp\Big\{ -{\sum_{i=1}^k\Vert\x_i-\bnuh\Vert_{\V_i^{-1}}^2 \over 2(\tau^2+\si^2)} - {S\over 2\si^2}\Big\}.
\end{align*}
Then, the hierarchical Bayes estimator of $\bmu_1$ relative to the quadratic loss (\ref{eqn:loss}) is written as
\begin{align}
\bmuh_1^{HB} =& E[\bmu_1/\si^2\mid \X_1,  \ldots, \X_k, S]/E[1/\si^2\mid \X_1,  \ldots, \X_k, S]\non\\
=&
\X_1 - {E[(\tau^2+\si^2)^{-1} \mid \X_1,  \ldots, \X_k, S] \over E[(\si^2)^{-1} \mid \X_1,  \ldots, \X_k, S] }(\X_1 - \bnuh).
\label{eqn:HB}
\end{align}
Making the transformation $\la=\si^2/(\tau^2+\si^2)$ and $\eta=1/\si^2$ with the Jacobian $|\partial(\tau^2, \si^2)/\partial (\la,\eta)|=1/(\la^2\eta^3)$ gives
$$
\pi(\la, \eta \mid \x_1, \ldots, \x_k,S)
\propto  \la^{p(k-1)/2+a -1}\eta^{\{n+p(k-1)\}/2-c -1}\exp\Big\{ - {\la\eta\over 2}\sum_{i=1}^k\Vert\x_i-\bnuh\Vert_{\V_i^{-1}}^2 -{\eta\over 2}S\Big\},
$$
where $0<\la<1$ and $\eta\geq L$.
Thus, we have
\begin{align*}
{E[(\tau^2+\si^2)^{-1} \mid \X_1,  \ldots, \X_k, S] \over E[(\si^2)^{-1} \mid \X_1,  \ldots, \X_k, S] }
=&{\int_0^1 \int_L^\infty \la^{p(k-1)/2+a }\eta^{\{n+p(k-1)\}/2-c }\exp\{ - {\eta S\over 2}(\la F+1)\}d\eta d\la
\over
\int_0^1 \int_L^\infty \la^{p(k-1)/2+a -1}\eta^{\{n+p(k-1)\}/2-c }\exp\{ - {\eta S\over 2}(\la F+1)\}d\eta d\la}
\\
=& {\phi^{HB}(F,S) / F},
\end{align*}
for
\begin{equation}
\phi^{HB}(F,S)=
{\int_0^F \int_{LS}^\infty x^{p(k-1)/2+a }v^{\{n+p(k-1)\}/2-c }\exp\{ - {v}(x+1)/2\}dv dx
\over
\int_0^F \int_{LS}^\infty x^{p(k-1)/2+a -1}v^{\{n+p(k-1)\}/2-c }\exp\{ - {v}(x+1)/2\}dv dx},
\\
\label{eqn:phiHB}
\end{equation}
where the transformations $x=F\la$ and $v=S\eta$ are used.

\medskip
It is noted that the hierarchical Bayes estimator belongs to the class (\ref{eqn:class1}), and the minimaxity can be shown by checking the conditions (a) and (b) in Theorem \ref{thm:1}.
By differentiating $\phi^{HB}(F,S)$ with respect to $F$ and $S$, the condition (a) can be easily verified.
The condition (a) implies that
\begin{align*}
\phi^{HB}(F,S)\leq \lim_{F\to\infty}\lim_{S\to 0} \phi^{HB}(F,S)
=&
{\int_0^\infty x^{p(k-1)/2+a }/ (1+x)^{\{n+p(k-1)\}/2+1-c } dx
\over
\int_0^\infty x^{p(k-1)/2+a -1}/(1+x)^{\{n+p(k-1)\}/2+1-c }dx}\\
=& {B(p(k-1)/2+a +1, n/2-a  -c ) \over B(p(k-1)/2+a , n/2-a -c +1)}
= {p(k-1)+2a  \over n- 2(a +c )},
\end{align*}
for the beta function $B(a,b)$ if $a>-p(k-1)/2$ and $a+c<n/2$.
Thus, the condition (b) is satisfied if 
\begin{align}
&a>-p(k-1)/2, a+c<n/2,\nonumber \\
&{p(k-1)+2a  \over n- 2(a +c )}\leq 2 [\tr\{(\V_1-\A)\Q\}/{\rm Ch}_{\rm max}\{(\V_1-\A)\Q\}-2]/(n+2),
\label{eqn:condHB}
\end{align}
which provides conditions on $a$ and $c$ for the minimaxity of the hierarchical Bayes estimator.

\subsection{Hierarchical empirical Bayes minimax estimator under the normal prior for $\bnu$}
\label{sec:3.3}

The empirical Bayes estimator $\bmuh_1^{EB}$ and the hierarchical Bayes estimator $\bmuh_1^{HB}$ are derived under the uniform prior distribution for $\bnu$.
Instead of the uniform prior, we here assume the normal prior distribution for $\bnu$, namely,
\begin{equation}
\begin{split}
\bmu_i\mid \bnu, \tau^2 \sim& \Nc_p(\bnu, \tau^2\V_i), \quad \text{for}\ i=1, \ldots, k,
\\
\bnu \mid \ga^2,  \sim& \Nc_p(0, \ga^2 \A),
\end{split}
\label{eqn:prior3}
\end{equation}
for $\ga>0$ and  $\A=\Big(\sum_{i=1}^k\V_i^{-1}\Big)^{-1}$.
Then the posterior distributions are given by
\begin{equation}
\begin{split}
\bmu_i \mid \X_i, \bnu, \tau^2, \si^2 \sim& \Nc_p\Big( \bmuh_i^*(\si^2, \tau^2,\bnu), (\si^{-2}+\tau^{-2})^{-1}\Big),\quad\quad  i=1, \ldots, k,
\\
\bnu \mid \X_1, \ldots, \X_k, \tau^2,\ga^2,\si^2 \sim& \Nc_p\Big( {\ga^2\over \ga^2 + \tau^2+\si^2}\bnuh, {\ga^2(\tau^2+\si^2)\over \ga^2 + \tau^2+\si^2}\A \Big),
\end{split}
\label{eqn:posterior3}
\end{equation}
where $\bmuh_i^*(\si^2, \tau^2,\bnu)$ is defined below (\ref{eqn:posterior1}).
Thus, the Bayes estimator is
\begin{align}
\bmuh_1^B(\si^2, \tau^2,\ga^2)=&
\X_1 - {\si^2\over \tau^2+\si^2}(\X_1-E[\bnu\mid \X_1, \ldots, \X_k])
\non\\
=&
\X_1 - {\si^2\over \tau^2+\si^2}(\X_1-\bnuh)-{\si^2\over \ga^2+\tau^2+\si^2}\bnuh.
\label{eqn:Bayes3}
\end{align}
Because the marginal density of $\X_1, \ldots, \X_k$ is
\begin{equation}
\begin{split}
f_\pi(\x_1,\ldots, \x_k \mid \tau^2, \ga^2, \si^2) \propto&
{1\over (\tau^2+\si^2)^{p(k-1)/2}}\exp\Big\{ -{\sum_{i=1}^k\Vert\x_i-\bnuh\Vert_{\V_i^{-1}}^2 \over 2(\tau^2+\si^2)} \Big\}
\\
&\times {1\over (\ga^2 + \tau^2+\si^2)^{p/2}} \exp\Big\{ - {\Vert\bnuh\Vert_{\A^{-1}}^2 \over 2(\ga^2 + \tau^2+\si^2)}\Big\},
\end{split}
\label{eqn:posterior3}
\end{equation}
we can estimate $\tau^2+\si^2$, $\ga^2+\tau^2+\si^2$ by $\sum_{i=1}^k \Vert \X_i-\bnuh\Vert_{\V_i^{-1}}^2/\{p(k-1)-2\}$ and $\Vert\bnuh\Vert_{\A^{-1}}^2/(p-2)$, respectively, from the marginal likelihood.
When $\si^2$ is estimated by $\sih^2=S/(n+2)$, the resulting hierarchical empirical Bayes estimator is
\begin{equation}
\bmuh_1^{HEB} = \X_1 - \min\Big( {a_0\over F}, 1\Big)(\X_1-\bnuh)
-\min\Big( {b_0\over G}, 1\Big)\bnuh,
\label{eqn:HEB}
\end{equation}
for $a_0=\{p(k-1)-2\}/(n+2)$ and $b_0=(p-2)/(n+2)$.
This estimator belongs to the class (\ref{eqn:class2}).
It follows from Theorem \ref{thm:2} that the hierarchical empirical Bayes estimator $\bmuh_1^{HEB}$ is minimax if $0<a_0\leq [\tr\{(\V_1-\A)\Q\}/{\rm Ch}_{\rm max}\{(\V_1-\A)\Q\}-2]/(n+2)$ and $0<b_0\leq\{\tr(\A\Q)/{\rm Ch}_{\rm max}(\A\Q)-2\}/(n+2)$.

\section{Extension to Estimation of Linear Combinations}
\label{sec:4}

We here extend the result in Section \ref{sec:2} to estimation of the linear combination of $\bmu_1, \ldots, \bmu_k$, namely,
$$
\bth=\sum_{i=1}^k d_i \bmu_i,
$$
where $d_1, \ldots, d_k$ are constants.
Based on (\ref{eqn:class1}), we consider the class of the estimators
\begin{equation}
\bthh(\phi) = \sum_{i=1}^k d_i \Big\{\X_i - {\phi(F, S)\over F}(\X_i-\bnuh)\Big\}.
\label{eqn:class3}
\end{equation}
When the estimator is evaluated in light of risk relative to the loss function $\Vert\bthh(\phi)-\bth\Vert^2_\Q/\si^2$, we obtain conditions for minimaxity of the estimators (\ref{eqn:class3}).

\begin{thm}
\label{thm:3}
Assume the conditions ${\rm Ch}_{\rm max}((\sum_{i=1}^k d_i^2\V_i-(\sum_{i=1}^k d_i)^2\A)\Q)\neq0$ and $\tr\{(\sum_{i=1}^k d_i^2\V_i-(\sum_{i=1}^k d_i)^2\A)\Q\}/{\rm Ch}_{\rm max}((\sum_{i=1}^k d_i^2\V_i-(\sum_{i=1}^k d_i)^2\A)\Q)>2$. Then, the estimator $\bthh(\phi)$ is minimax  if $\phi(F,S)$ satisfies the following conditions:

{\rm (a)}\ $\phi(F,S)$ is non-decreasing in $F$ and non-increasing in $S$.

{\rm (b)}\ $\phi(F,S)$ satisfies the inequality
\begin{equation}
0<\phi(F,S)\leq {2\over n+2} \Big[ {\tr\{(\sum_{i=1}^k d_i^2\V_i-(\sum_{i=1}^k d_i)^2\A)\Q\}\over {\rm Ch}_{\rm max}((\sum_{i=1}^k d_i^2\V_i-(\sum_{i=1}^k d_i)^2\A)\Q)} -2 \Big].
\label{eqn:thm3b}
\end{equation}
\end{thm}

{\bf Proof}.\ \ 
The risk function of the estimator $\bthh(\phi)$ is 
\begin{align*}
R(\bom, \bthh(\phi))=&
\sum_{i=1}^k {d_i^2 \over \si^2}E\Big[ \Big\Vert \X_i - \bmu_i - {\phi(F, S)\over F}(\X_i-\bnuh) \Big\Vert_\Q^2\Big]\\
&+
\sum_{i\not= j} {d_id_j\over \si^2} 
E\Big[ \Big\{ \X_i - \bmu_i - {\phi(F, S)\over F}(\X_i-\bnuh) \Big\}^\top \Q\Big\{ \X_j - \bmu_j - {\phi(F, S)\over F}(\X_j-\bnuh) \Big\}\Big]
\\
=& K_1+K_2. \quad \text{ (say)}
\end{align*}
Concerning $K_2$, the same arguments as in (\ref{eqn:I2}) and (\ref{eqn:I3}) are used to get
\begin{align*}
{1\over \si^2}E\Big[ (\X_i-\bmu_i)^\top\Q(\X_j-\bnuh) {\phi\over F}\Big]
=& E\Big[\tr\Big[\V_i\Q\bnabla_i \Big\{(\X_j-\bnuh)^\top {\phi(F,S)\over F}\Big\}\Big]\Big]
\non\\
=&
E\Big[ - \tr(\A\Q){\phi \over F} - 2 B_{ij}{\phi \over F} + 2 B_{ij}\phi_F\Big], 
\end{align*}
and
\begin{align*}
B_{ij} E\Big[{S\over \si^2}{\phi^2\over F}\Big] 
=&
B_{ij} E\Big[ (n+2) {\phi^2\over F} - 4\phi\phi_F + 4S {\phi\phi_S\over F} \Big],
\end{align*} 
where 
\begin{equation}
B_{ij}= (\X_i-\bnuh)^\top\Q(\X_j-\bnuh)/\sum_{a=1}^k (\X_a-\bnuh)^\top\V_a^{-1}(\X_a-\bnu).
\label{eqn:Bij}
\end{equation}
Thus, $K_2$ is written as
\begin{equation}
K_2=\sum_{i\not= j} d_id_j\Big[ 2\tr(\A\Q){\phi\over F} 
+ B_{ij} \Big\{4 {\phi\over F} - 4 \phi_F + (n+2){\phi^2\over F}- 4 \phi\phi_F+4S{\phi\phi_S\over F}\Big\}\Big].
\label{eqn:K2}
\end{equation}
Concerning $K_1$, from (\ref{eqn:p2}), it follows that
\begin{align}
K_1=&
\sum_{i=1}^k d_i^2 E\Big[ \tr(\V_i\Q) - 2 \tr\{(\V_i-\A)\Q\}{\phi\over F} \non\\
&+ B_{ii} \Big\{ (n+2){\phi^2\over F} + 4 {\phi\over F} 
-4\phi_F-4\phi\phi_F+4S {\phi\phi_S\over F}\Big\} \Big].
\label{eqn:K1}
\end{align}
It is here observed that
$$
\sum_{i=1}^k d_i^2 B_{ii} + \sum_{i\not= j}d_id_j B_{ij}
=B(\X),
$$
where $\X=(\X_1, \ldots, \X_k)$ and
\begin{equation}
B(\x_1, \ldots, \x_k) = {\{\sum_{i=1}^k d_i(\x_i-\bnuh)\}^\top\Q\{\sum_{i=1}^k d_i(\x_i-\bnuh)\} \over
\sum_{j=1}^k (\x_j-\bnuh)^\top\V_j^{-1}(\x_j-\bnu)}.
\label{eqn:B2}
\end{equation}
Combining (\ref{eqn:K2}) and (\ref{eqn:K1}), one gets
\begin{align*}
R(\bom, \bthh(\phi))=&
\sum_{i=1}^k d_i^2 \tr(\V_i\Q) 
\\
&+ E\Big[ {\phi\over F} \Big\{- 2 \sum_{i=1}^k\tr(d_i^2\V_i\Q) + 2 \Big(\sum_{i=1}^k d_i\Big)^2 \tr(\A\Q) + 4 B(\X) + (n+2) \phi B(\X) \Big\}
\\
&+ B(\X)\Big\{ -4\phi_F-4\phi\phi_F+4S {\phi\phi_S\over F}\Big\} \Big],
\end{align*}
which is smaller than $R(\bom, \sum_{i=1}^k d_i\X_i)$ under the conditions (a) and (b) in Theorem \ref{thm:3}, because $B(\x)\leq {\rm Ch}_{\rm max}((\sum_{i=1}^k d_i^2\V_i-(\sum_{i=1}^k d_i)^2\A)\Q)$ from Lemma \ref{lem:in2}.
Hence, the proof of Theorem \ref{thm:3} is complete.
\hfill$\Box$

\begin{lem}
\label{lem:in2}
For $B(\x_1, \ldots, \x_k)$ given in $(\ref{eqn:B})$, it holds that
$$
B(\x_1, \ldots, \x_k)\leq
{\rm Ch}_{\rm max}\Big( \Big(\sum_{i=1}^k d_i^2\V_i-(\sum_{i=1}^k d_i)^2\A\Big)\Q\Big).
$$
\end{lem}

{\bf Proof}.\ \ 
Let $\y_i=\x_i-\bnuh$ and $\x=(\x_1, \ldots, \x_k)$.
Then, it is noted that $\sum_{i=1}^k \V_i^{-1}\y_i=\sum_{i=1}^k \V_i^{-1}(\x_i-\bnuh)=\zero$, which means that
$$
\W \begin{pmatrix}\y_1\\ \vdots\\ \y_k \end{pmatrix}=\begin{pmatrix}\y_1\\ \vdots\\ \y_k \end{pmatrix},
$$
for
$$
\W = {\rm block\ diag}(\I, \ldots, \I) - \begin{pmatrix}\A\\ \vdots\\ \A\end{pmatrix}(\V_1^{-1}, \ldots, \V_k^{-1}).
$$ 
Then, the numerator and the denominator of $B(\x)$ are
\begin{align*}
&\{\sum_{i=1}^k d_i(\x_i-\bnuh)\}^\top\Q\{\sum_{i=1}^k d_i(\x_i-\bnuh)\} 
= (\y_1^\top, \ldots, \y_k^\top) \begin{pmatrix}d_1\I\\ \vdots \\ d_k\I\end{pmatrix}\Q (d_1\I, \ldots, d_k \I)\begin{pmatrix}\y_1\\ \vdots\\ \y_k\end{pmatrix}
\\
&\quad
= (\y_1^\top, \ldots, \y_k^\top) \W^\top \begin{pmatrix}d_1\I\\ \vdots \\ d_k\I\end{pmatrix}\Q(d_1\I, \ldots, d_k \I)\W \begin{pmatrix}\y_1\\ \vdots\\ \y_k\end{pmatrix},
\\
&\sum_{j=1}^k (\x_j-\bnuh)^\top\V_j^{-1}(\x_j-\bnu)
=  (\y_1^\top, \ldots, \y_k^\top) {\rm block\ diag}(\V_1^{-1}, \ldots, \V_k^{-1})\begin{pmatrix}\y_1\\ \vdots\\ \y_k\end{pmatrix}.
\end{align*}
Thus, we get an upper bound given by
\begin{align*}
B(\x) \leq& {\rm Ch}_{\rm max}\Big( {\rm block\ diag}(\V_1, \ldots, \V_k) \W^\top \begin{pmatrix}d_1\I\\ \vdots \\ d_k\I\end{pmatrix}\Q(d_1\I, \ldots, d_k \I) \W \Big)
\\
=& {\rm Ch}_{\rm max}\Big( \Q(d_1\I, \ldots, d_k \I) \W {\rm block\ diag}(\V_1, \ldots, \V_k) \W^\top \begin{pmatrix}d_1\I\\ \vdots \\ d_k\I\end{pmatrix} \Big)
\\
=&
{\rm Ch}_{\rm max}\Big( \Big(\sum_{i=1}^k d_i^2\V_i-(\sum_{i=1}^k d_i)^2\A\Big)\Q\Big),
\end{align*}
which shows Lemma \ref{lem:in2}.
\hfill$\Box$

\bigskip
We here give some examples of the condition (\ref{eqn:thm3b}) in specific cases.
For example, in the case of $k=2$, it is observed that $\tr\{(\sum_{i=1}^2 d_i^2\V_i-(\sum_{i=1}^2 d_i)^2\A)\Q\}=\tr\{ (d_1\V_1-d_2\V_2)(\V_1+\V_2)^{-1}(d_1\V_1-d_2\V_2)\Q\}=\tr(\H)$, where
$$
\H=(\V_1+\V_2)^{-1/2}(d_1\V_1-d_2\V_2)\Q (d_1\V_1-d_2\V_2)(\V_1+\V_2)^{-1/2}.
$$
Also, it can be seen that ${\rm Ch}_{\rm max}((\sum_{i=1}^k d_i^2\V_i-(\sum_{i=1}^k d_i)^2\A)\Q)={\rm Ch}_{\rm max}(\H)$.
Hence, the condition (b) in Theorem \ref{thm:3} is expressed as
$$
0<\phi(F,S)\leq {2\over n+2} \Big[ {\tr(\H)\over {\rm Ch}_{\rm max}(\H)} -2 \Big].
$$

For example, in the case that $\V_1=\cdots =\V_k=\Q=\I$, we have
$$
\tr\{(\sum_{i=1}^k d_i^2\V_i-(\sum_{i=1}^k d_i)^2\A)\Q\}
= p \sum_{i=1}^k (d_i-\ddo)^2,
$$
for $\ddo=k^{-1}\sum_{i=1}^k d_i$.
Similarly, ${\rm Ch}_{\rm max}((\sum_{i=1}^k d_i^2\V_i-(\sum_{i=1}^k d_i)^2\A)\Q)=\sum_{i=1}^k (d_i-\ddo)^2$, which implies that the condition (b) is expressed as
$$
0<\phi(F,S)\leq 2(p-2)/(n+2).
$$

\section{Simulation Studies}
\label{sec:sim}

We investigate the numerical performances of the risk functions of the preliminary-test estimator and several empirical and hierarchical Bayes estimators through simulation.
We employ the quadratic loss function $L(\bde_1, \bmu_1, \si^2)$ in (\ref{eqn:loss}) for $\Q=\V_1^{-1}$.

\medskip
The estimators which we compare are the following five:

\medskip
PT: the preliminary-test estimator given in (\ref{eqn:PT})
$$
\bmuh_1^{PT}= 
\left\{\begin{array}{ll} \X_1 & \ \text{if}\ F>(p(k-1)/n)F_{p(k-1), n, \al},\\
\bnuh & \ \text{otherwise},
\end{array}\right.
$$

JS: the James-Stein estimator
$$
\bmuh_1^{JS}=\X_1 - {p-2 \over n+2}{S\over \Vert\X_1\Vert^2_{\V_1^{-1}}}\X_1,
$$

EB: the empirical Bayes estimator given in (\ref{eqn:EB1}) 
$$
\bmuh_1^{EB} = \X_1 - \min\Big( {a_0 \over F}, 1\Big)(\X_1-\bnuh),
$$
for $a_0= [\tr\{(\V_1-\A)\V_1^{-1}\}/{\rm Ch}_{\rm max}\{(\V_1-\A)\V_1^{-1}\}-2]/(n+2)$ (one can see that this constant choice is optimal with respect to the upper bound of the risk difference), 

HB: the hierarchical Bayes estimator given in (\ref{eqn:HB}) and (\ref{eqn:phiHB}), 
$$
\bmuh_1^{HB} =\X_1 - {\phi^{HB}(F,S) \over F}(\X_1 - \bnuh),
$$

HEB: the hierarchical empirical Bayes estimator given in (\ref{eqn:HEB}) 
$$
\bmuh_1^{HEB} = \X_1 - \min\Big( {a_0\over F}, 1\Big)(\X_1-\bnuh)
-\min\Big( {b_0\over G}, 1\Big)\bnuh,
$$
for $a_0= [\tr\{(\V_1-\A)\V_1^{-1}\}/{\rm Ch}_{\rm max}\{(\V_1-\A)\V_1^{-1}\}-2]/\{2(n+2)\}$ and $b_0=\{\tr(\A\V_1^{-1})/{\rm Ch}_{\rm max}(\A\V_1^{-1})-2\}/\{2(n+2)\}$ (these constants are also optimal choices).

\medskip
It is noted that the James-Stein estimator does not use $\X_2, \ldots, \X_k$, but is minimax.
Concerning the hierarchical Bayes estimator $\bmuh_1^{HB}$, the constants $c$ and $L$ are $c=1$ and $L=0$, and $a$ is the solution of the equation
$$
{p(k-1)+2a  \over n- 2(a +1 )}(n+2)= \tr\{(\V_1-\A)\V_1^{-1}\}/{\rm Ch}_{\rm max}\{(\V_1-\A)\V_1^{-1}\}-2,
$$
which guarantees the minimaxity from the condition (\ref{eqn:condHB}).

\medskip
In this simulation, we generate random numbers of $\X_1, \ldots, \X_k$ and $S$ based on the model (\ref{eqn:model}) for $p=k=5$, $n=20$, $\si^2=2$ and $\V_i=(0.1\times i)\I_p$, $i=1, \ldots, k$.
For the mean vectors $\bmu_i$, we treat the 12 cases: 
\begin{align*}(\bmu_1,& \ldots, \bmu_5)
\\
=&(\zero, \zero, \zero, \zero, \zero), (1\j_5, 1\j_5, 1\j_5, 1\j_5, 1\j_5), (2\j_5, 2\j_5, 2\j_5, 2\j_5, 2\j_5), (3\j_5, 3\j_5, 3\j_5, 3\j_5, 3\j_5),\\ 
&(-0.4\j_5, -0.2\j_5, \zero, 0.2\j_5, 0.4\j_5), (2\j_5, -0.5\j_5, -0.5\j_5, -0.5\j_5, -0.5\j_5),  \\
&(4\j_5, -1\j_5, -1\j_5, -1\j_5, -1\j_5), \\
& (1.2\j_5, 1.4\j_5, 1.6\j_5, 1.8\j_5, 2\j_5),(\zero, 2\j_5, 2\j_5, 2\j_5, 2\j_5), (\zero, 4\j_5, 4\j_5, 4\j_5, 4\j_5), (2\j_5, \zero, \zero, \zero, \zero, \zero), \\ 
\end{align*}
where $\j_p=(1, \ldots, 1)^\top\in \Re^p$.
The first four are the cases of equal means, the next three are the cases with $\sum_{i=1}^5\bmu_i=\zero$ and the last four are various unbalanced cases.

\medskip
For each estimator $\bmuh_1$, based on 5,000 replication of simulation, we obtain an approximated value of the risk function $R(\bom, \bmuh_1)=E[L(\bmuh_1, \bmu_1, \si^2)]$.
Table \ref{table:risk1} reports the percentage relative improvement in average loss (PRIAL) of each estimator $\bmuh_1$ over $\X_1$, defined by
$$
{\rm PRIAL} = 100\times\{ R(\bom, \X_1) - R(\bom, \bmuh_1)\}/R(\bom, \X_1).
$$

\small
\begin{table}[!thb]
\caption{Values of PRIAL of estimators PT, JS, EB, HB and HEB}
\begin{center}
$
{\renewcommand\arraystretch{1.1}\small
\begin{array}{c@{\hspace{5mm}}
              r@{\hspace{2mm}}
              r@{\hspace{2mm}}
              r@{\hspace{2mm}}
              r@{\hspace{2mm}}
              r@{\hspace{2mm}}
              r@{\hspace{2mm}}
              r
             }
\text{$(\bmu_1, \ldots, \bmu_5)$} &\text{PT}&\text{JS}&\text{EB}&\text{HB}&\text{HEB}\\
\hline
\text{$(0,0,0,0,0)\otimes \j_5$} 
&52.15317
&53.97469
&14.66425
&14.57437
&26.38606
\\
\text{$(1, 1, 1, 1, 1)\otimes \j_5$} 
&52.15317
&12.89115
&14.66425
&14.57437
&9.891098
\\
\text{$(2, 2, 2, 2, 2)\otimes \j_5$} 
&52.15317
&4.066823
&14.66425
&14.57437
&8.516356
 \\
 \text{$(3, 3, 3, 3, 3)\otimes \j_5$} 
 &52.15317
 &2.268442
 &14.66425
 &14.57437
 &8.249028
 \\
\hline
\text{$(-0.4,-0.2,0,0.2,0.4)\otimes \j_5$} 
&37.34717
&37.20396
&13.01833
&12.97352
&22.64692
\\
\text{$(2,-0.5,-0.5,-0.5,-0.5)\otimes \j_5$} 
&-56.8291
&4.066823
&3.213333
&3.213459
&6.031053
\\
 \text{$(4,-1,-1,-1,-1)\otimes \j_5$} 
 &0.7375904
 &1.620614
 &1.358956
 &1.358821
 &2.098222
 \\
\hline
\text{$(1.2,1.4,1.6,1.8,2)\otimes \j_5$} 
&37.34717
&9.463467
&13.01833
&12.97352
&8.397694
 \\
\text{$(0.2, 2, 2, 2, 2)\otimes \j_5$}
&-98.94453
&49.73141
&4.947591
&4.949466
&5.183324 
 \\
\text{$(0.4, 4,4,4,4)\otimes \j_5$} 
&-2.492994
&37.56052
&1.805795
&1.80584
&2.071347
 \\
\text{$(2,0,0,0,0)\otimes \j_5$} 
&-94.45962
&4.066823
&4.439434
&4.440511
&4.479298
\\
\end{array}
}
$
\end{center}
\label{table:risk1}
\end{table}
\normalsize

It is revealed from Table \ref{table:risk1} that the performance of the preliminary test estimator PT strongly depends on the setup of parameters, namely it is good under the hypothesis of equal means, but not good for parameters close to the hypothesis.
The James-Stein estimator JS is good for small $\bmu_1$, but not good for large $\bmu_1$.
The empirical Bayes estimator EB and the hierarchical Bayes estimator HB perform similarly and they remain good even for large $\bmu_1$ as long as $\bmu_1=\cdots\bmu_5$.
The performance of HEB depends on parameters and is good for smaller means.
For means with $\sum_{i=1}^5 \bmu_i=\zero$, HEB is better than EB and HB.
Thus, EB, HB and HEB are used as an alternative to PT.

\section{Concluding Remarks}
\label{sec:remark}

An interesting query is to find an admissible and minimax estimator of $\bmu_1$.
In the framework of simultaneous estimation of $(\bmu_1, \ldots, \bmu_k)$, Imai, $\et$ (2017) demonstrated that all the estimators within the class (\ref{eqn:class1}) are improved on by an estimator belonging to the class (\ref{eqn:class2}), which means that the hierarchical Bayes estimators against the uniform prior of $\bnu$ is inadmissible.
However, we could not show the same story in the single estimation of $\bmu_1$.
Because the hierarchical Bayes estimator $\bmuh_1^{HB}$ is derived under the uniform prior for $\bnu$, it could be supposedly inadmissible. 
Thus, whether it is admissible or not is an open question.

\medskip
An approach to the admissible and minimax estimation is the proper prior distribution
\begin{equation}
\begin{split}
\pi(\tau^2 \mid \si^2) \propto& \Big({\si^2\over \tau^2+\si^2}\Big)^{a +1},\\
\pi(\ga^2 \mid \tau^2, \si^2) \propto& \Big({\si^2\over \ga^2 + \tau^2+\si^2}\Big)^{b +1},\\
\pi(\si^2) \propto& (\si^2)^{c -3},\quad \text{for}\ \si^2\leq 1/L ,
\end{split}
\label{eqn:prior4}
\end{equation}
where $a $, $b$ and $c $ are constants and $L $ is a positive constant.
As seen from Imai, $\et$ (2017), the resulting Bayes estimator has the form
\begin{align}
\bmuh_1^{FB} = 
\X_1 - {\phi^{FB}(F,G,S) \over F}(\X_1 - \bnuh) - {\psi^{FB}(F,G, S) \over G}\bnuh,
\label{eqn:FB}
\end{align}
where $\phi^{FB}(F,G,S)$ and $\psi^{FB}(F,G,S)$ are functions of $F$, $G$ and $S$.
Unfortunately, we could not establish the minimaxity for this type of estimators, which is another interesting query.

\medskip
In this paper, we investigated the minimaxity of Bayesian alternatives to the preliminary test estimator.
Beyond this framework, we consider the estimation of $\bmu_1$ based on $\X_1$ and $S$ without $\X_2, \ldots, \X_k$.
Using the argument as in Strawderman (1973), we can derive an admissible and minimax estimator of $\bmu_1$.
Taking the hypothesis $H_0 : \bmu_1=\cdots=\bmu_k$ into account, we can use the same argument as in Strawderman (1973) under the assumption that $\X_2, \ldots, \X_k$ are given and fixed.
Namely, we consider the prior distributions (\ref{eqn:prior1}) and (\ref{eqn:prior2}) where $\bnu$ is replaced with $\bnu^*=(\sum_{i=2}^k \V_i^{-1})^{-1}\sum_{i=2}^k\V_i^{-1}\X_i$.
Then, the Bayes estimator is 
$$
\X_1 - {\phi^S(\Vert\X_1-\bnu^*\Vert^2, S)\over \Vert\X_1-\bnu^*\Vert^2/ S}(\X_1-\bnu^*),
$$
where $\phi^S$ is a function derived from the Bayes estimator against the Strawderman type prior.
This shrinks $\X_1$ towards $\bnu^*$ and is minimax under some condition on $\phi^S$.
Thus, it is admissible and minimax in the framework that $\X_2, \ldots, \X_k$ are given and fixed.

\section*{Acknowledgments}
We would like to thank the Editor, the Associate Editor and the reviewer for valuable comments and helpful suggestions which led to an improved version of this paper.
Research of the second author was supported in part by Grant-in-Aid for Scientific Research  (15H01943 and 26330036) from Japan Society for the Promotion of Science.

\end{document}